\documentclass[11pt]{article}
\input{epsf.sty}
\usepackage{hyperref}
\usepackage{pstricks}
\usepackage{amsmath, amssymb} 
\newtheorem {thm}{Theorem}[section]
\newtheorem {prop}[thm]{Proposition}
\newtheorem {lem}[thm]{Lemma}

\newtheorem {defn}[thm]{Definition}

\def\Cox{\hfill \Box}

\def\P{{\Bbb P}}
\def\Q{{\Bbb Q}}
\def\E{{\Bbb E}}

\def\0{{\bf 0}}

\def\a{\alpha}

\def\b{\beta}

\def\d{\delta}

\def\e{\varepsilon}

\def\phi{\varphi}

\def\g{\gamma}

\def\l{\lambda}

\def\k{\kappa}

\def\s{\sigma}

\def\x{\xi}

\def\T{\T}

\def\LL{{\cal L}}

\begin{document}

\title{A symmetric entropy bound on 
the non-reconstruction regime of Markov chains on Galton-Watson trees } 

\author{
Marco Formentin
\footnote{
Universit\`a degli Studi di Padova,
Dipartimento di Matematica Pura ed Applicata,
Via Trieste, 63 - 35121 Padova, 
Italia
\texttt{formen@math.unipd.it}}
\, and
Christof K\"ulske
\footnote{
University of Groningen, 
Department of Mathematics and Computing Sciences, 
Nijenborgh 9,   
9747 AC Groningen, 
The Netherlands
\texttt{kuelske@math.rug.nl},
\texttt{ http://www.math.rug.nl/$\sim$kuelske/ }}
}

\maketitle

\begin{abstract}
We give a criterion of the form $\Q(d) c(M)<1$
for the non-reconstructability of tree-indexed $q$-state Markov chains 
obtained by broadcasting a signal from the root with a given transition matrix $M$. 
Here $c(M)$ is an explicit function, which is convex over the set of $M$'s with a given invariant distribution,  that is defined in terms of  a $q-1$-dimensional 
variational problem over symmetric entropies. Further 
$\Q(d)$ is the expected number of offspring on the Galton-Watson tree. 

This result is equivalent to proving the extremality of the free boundary 
condition Gibbs measure within the corresponding Gibbs-simplex.  

Our theorem holds for possibly non-reversible $M$ and its proof 
is based on a general  recursion formula 
for expectations of a symmetrized relative entropy function, 
which invites their use as a Lyapunov function. 
In the case of the Potts model, the present theorem reproduces 
earlier results of the authors, with a simplified proof, 
in the case of the symmetric Ising model (where the argument becomes similar 
to the approach of Pemantle and Peres) the method produces 
the correct reconstruction threshold), 
in the case of the (strongly) asymmetric Ising model where the Kesten-Stigum bound 
is known to be not sharp the method provides 
improved numerical bounds.

\end{abstract}

\smallskip
\noindent {\bf AMS 2000 subject classification:} 60K35, 82B20, 82B44 
\bigskip 

{\em Keywords:} Broadcasting on trees, 
Gibbs measures, random tree, Galton-Watson tree,
 reconstruction problem, free boundary condition.

\section{Introduction} \label{sect:intro}

The problem of reconstruction of Markov chains on $d$-ary trees has enjoyed 
much interest in recent years. There are multiple reasons for this, one of them 
being that it is a topic where people from information theory, {researchers in mathematical statistical mechanics,} pure probabilists, 
and people from the theoretical physics 
side of statistical mechanics can meet and make contributions. 

Indeed, starting with the symmetric Ising channel for which the reconstruction 
threshold was settled in  \cite{BlRuZa95, Io96a,Io96b}, using different methods 
and increased generality w.r.t. the underlying tree, there have been publications 
by a.o. Borgs, Chayes, Janson, Mossel, Peres from the mathematics side
\cite{Mo01, MoPe03, JaMo04, Bo06},  
deriving upper and lower bounds on reconstruction threshold for certain models 
of finite state tree-indexed Markov chains. 
From the theoretical physics side let us highlight \cite{MeMo06} on trees 
(see also \cite{GeMo07} on graphs) which contains a 
discussion of the potential 
relevance of  the reconstruction problem also to the glass problem. That paper  
also provides numerical values for the Potts model on the basis of 
extensive simulation results. The Potts model is interesting because, unlike 
the Ising model,  the true reconstruction threshold behaves (respectively is expected to behave) 
non-trivially as a function of the degree $d$ of the underlying $d$-ary tree and the number of states $q$. 
For a discussion of this see the conjectures in \cite{MeMo06}, the rigorous 
bounds in \cite{FoKu08}, 
and in particular the proof in \cite{Sl08} 
showing that the Kesten-Stigum bound is not sharp if $q\geq 5$, and sharp if $q=3$, for large enough $d$.  
{ We refer to \cite{BhMa09} for a general computational method to obtain non trivial rigorous 
bounds for reconstruction on trees.}

Now, our treatment is motivated in the generality of its setup 
by the questions raised and type of results given 
in \cite{Ma04}, and technically somewhat inspired by \cite{PePe06, FoKu08}. 
Indeed, for the Potts model 
the present paper reproduces the result of \cite{FoKu08} (where moreover  
 we provide numerical estimates on the reconstruction inverse temperature, also in the small $q$, small 
 $d$ regime.)   
However, in the present paper the focus is on generality, that is the universality of the type of 
estimate, and the structural clarity of the proof.  
It should be clear that the condition we provide can be easily implemented 
in any given model to produce numerical estimates on reconstruction thresholds.  

The remainder of the paper is organized as follows.  
Section  \ref{s2} contains the definition of the model and the statement of the theorem.   
Section  \ref{s3} contains the proof.

\bigskip 
\bigskip

\section{Model and result}\label{s2}

Consider an infinite rooted tree $T$ having no leaves.  
For $v,w\in T$ we write $v \rightarrow w$, if $w$ is the child of $v$,  and we denote by 
$|v|$ the distance of a vertex $v$ to the root. 
We write $T^N$ for the subtree of all vertices with distance $\leq N$ to the root. 


To each vertex $v$ there is associated a (spin-) variable 
$\sigma(v)$ taking values in a finite space which, 
without loss of generality, will be denoted by $\{1,2,\dots,q\}$.
Our model will be defined in terms of the stochastic matrix with non-zero entries 
\begin{equation}\begin{split}\label{transegeneral}
M=(M(v,w))_{1\leq v,w\leq q}\:.
\end{split}
\end{equation} 
By the Perron-Frobenius theorem there is a unique  
single-site measure $\alpha=(\alpha(j))_{j=1,\dots,q}$ which is invariant 
under the application of the transition matrix $M$,  meaning that 
$\sum_{i=1}^q\alpha(i)M(i,j)=\alpha(j)$. 

The object of our study is the  corresponding 
{\em tree-indexed Markov chain in equilibrium}. This is the probability 
distribution $\P$ on $\{1,\dots,q\}^{T}$ whose restrictions $\P_{T^N}$ to the 
state spaces of finite trees 
$\{1,\dots,q\}^{T^N}$ are given by 
\begin{equation}\begin{split}\label{2}
&\P_{T^N}(\sigma_{T^N})=\alpha(\sigma(0))\prod_{{v,w:}\atop{ v\rightarrow w}} M(\sigma(v),\sigma(w))\:.\cr
\end{split}
\end{equation} 
The notion equilibrium refers to the fact 
that all single-site marginals are given by the invariant measure $\alpha$. 

A probability measure $\mu$ on $\{1,2,\dots,q\}^T$ is 
called a {\em Gibbs measure} if it  has the same finite-volume conditional probabilities 
as $\P$ has.  
This means that,  for all finite subsets $V\subset T$, we have for all $N$ sufficiently large   
\begin{equation}\begin{split}\label{2}
&\mu(\s_V | \s_{V^c})= \P_{T^N}(\s_{V}|\s_{\partial V})
\end{split}
\end{equation} 
$\mu$-almost surely.
The Gibbs measures, being defined in terms of a linear equation,  
form a simplex, and we would like to understand its structure, 
and exhibit its extremal elements \cite{Ge88}. Multiple Gibbs measures (phase transitions) 
may occur if the loss of memory in the transition described by $M$ 
is small enough compared to the proliferation of offspring along the tree $T$. 
Uniqueness of the Gibbs measure trivially implies extremality of the measure $\P$, 
but interestingly the converse is not true. Parametrizing $M$ by a temperature-like 
parameter may lead to two different transition temperatures, one where $\P$ 
becomes extremal and one where the Gibbs measure becomes unique.  
Broadly speaking, 
statistical mechanics models with two transition temperatures are peculiar to trees 
(and more generally to models indexed by non-amenable graphs \cite{Ly00}). This is 
one of the reasons for our interest in models on trees.

Now, our present aim is to provide a general criterion, depending 
on the model only in a local (finite-dimensional) way, which implies 
the extremality of $\P$, and which works also in regimes of non-uniqueness. 
People with statistical mechanics background may think 
of it as an analogy to Dobrushin's theorem saying that $c(\g)<1$ 
implies the uniqueness of the Gibbs measure of a local specification $\g$ 
where $c(\g)$ is determined in terms of local (single-site) quantities. 

In fact,  Martinelli et al. \cite{Ma04} (see Theorem 9.3., see also 
Theorem 9.3.' and Theorem 9.3'') give such a theorem. 
Their criterion for non-reconstruction of Markov chains on $d$-ary trees 
has the form $d \l_2 \k<1$ where $\k$ is  
the Dobrushin constant \cite{Ge88} of the system of 
conditional probabilities described by $\P$. Further $ \l_2 $ is the 
second eigenvalue of $M$. Our theorem takes a different form. 
Now, to formulate our result we need the following notation. 

We write for the simplex of length-$q$ probability vectors
\begin{equation}
\begin{split}\label{PP}
P=\{(p(i))_{i=1,\dots,q}, p(i) \geq 0\,\, \forall i,\,\, \sum_{i=1}^q p(i) =1\}
\end{split}
\end{equation}
and we denote 
the relative entropy between probability vectors $p,\a\in P$ by  
 $S(p|\a)=\sum_{i=1}^qp(i)\log\frac{p(i)}{\a(i)}$. We introduce the 
 {\em symmetrized entropy} between $p$ and $\a$ and write 
\begin{equation}
\begin{split}\label{defL}
L(p)=S(p|\a)+S(\a|p)=
{\sum_{i=1}^q \left(p(i)-\a(i)\right)\log\frac{p(i)}{\a(i)}}\:.
\end{split}
\end{equation}
While the symmetrized entropy is not a metric (since the triangle inequality fails) 
it serves us as a``distance" to the invariant measure $\a$. 

Let us define the constant, depending solely on the transition matrix $M$, 
in terms of the following supremum over probability vectors 
\begin{equation}\begin{split}\label{C1}
c(M)&=\sup_{p\in P}\frac{L(p M^\text{rev})}{L(p)},\cr
\end{split}
\end{equation}
where $M^\text{rev}(i,j)=\frac{\a(j)M(j,i)}{\a(i)}$ is the transition matrix of the reversed chain. 
Note that numerator and denominator vanish when 
we take for $p$ the invariant distribution $\a$. 
Consider a Galton-Watson tree with i.i.d. offspring distribution concentrated on $\{1,2,\dots\}$ 
and denote the corresponding expected number of offspring by $\Q (d)$.

Here is our main result. 

\begin{thm}\label{Theorem} If $\Q (d) c(M)<1$ then the tree-indexed Markov chain $\P$ 
on the Galton-Watson tree $T$ is extremal for $\Q$-almost every tree $T$. 
Equivalently, in information theoretic language, there is no reconstruction. 
\end{thm}

{\bf Remark 1.} The computation of the constant $c(M)$ for 
a given transition matrix $M$ is a simple numerical task. 
Note that fast mixing of the Markov chain corresponds to small $c(M)$. In this 
sense $c(M)$ is an effective quantity depending on the interaction $M$
that parallels the definition the Dobrushin constant $c_{D}(\g)$ in the theory of Gibbs measures measuring the degree of dependence in a local specification.  
While the latter depends on the structure of the interaction graph, 
this is {\em not} the case for $c(M)$. 

{\bf Remark 2.} Non-uniqueness of the Gibbs measures corresponds to the existence of boundary conditions 
which will cause the corresponding finite-volume conditional probabilities to converge 
to different limits. 
Extremality of the measure $\P$ means that conditioning the measure $\P$ 
to acquire a configuration $\xi$ at a distance larger than $N$ will cease 
to have an influence on the state at the root if $\xi$ is chosen according 
to the measure $\P$ itself and $N$ is tending to infinity. In the language 
of information theory this is called non-reconstructability (of the state at the 
origin on the basis of noisy observations far away). 

{\bf Remark 3. (On irreversibility.)} 
If $M$ is any transition matrix reversible for the equidistribution and, for a permutation 
$\pi$ of the numbers from $1$ to $q$,  we define $M_{\pi}(i,j)=M(i,\pi^{-1}j)$,  
then $c(M)=c(M_{\pi})$ for all permutations $\pi$.  This is seen by a simple computation.  
We can say that an irreversibility in the Markov chain which is caused by a deterministic 
stepwise renaming of labels (by $\pi$) is not seen in the constant.

{\bf Remark 4. (On Convexity.)} For all fixed probability vectors $\a$ 
the function $M\mapsto c(M)$ is 
convex on the set of transition matrices which 
have $\a$ as their invariant distribution, i.e. $\a M=\a$. 

This is a consequence of the fact that, for $M_1, M_2$ 
with $\a M_1=\a$, $\a M_2=\a$ we have 
that $(\l M_1+(1-\l)M_2)^\text{rev}= \l M^\text{rev}_1+(1-\l)M^\text{rev}_2$ and that the relative entropy is convex in the first 
and second argument. 

This implies that, for each $\a$, and fixed degree $d$, 
$\{M, \a M=\a; d c(M)<1\}$, for which  the criterion ensures 
non-reconstruction, is convex. 
\bigskip

We conclude this introduction with the discussion of two main types of test-examples. 

\bigskip

{\bf Example 1 (Symmetric Potts and Ising model.)} The Potts model with $q$ states 
at inverse temperature $\b$ is defined by the transition matrix 
\begin{equation}\begin{split}\label{transe}
M_{\b}=\frac{1}{e^{2\b}+ q-1}\begin{pmatrix}e^{2 \b} & 1& 1 & \dots & 1\cr 
1 & e^{2 \b} & 1 & \dots & 1\cr
1 & 1& e^{2 \b}  & \dots & 1\cr
&& \dots 
\end{pmatrix}
\end{split}.
\end{equation} 
This Markov chain is reversible for the equidistribution. In the 
case $q=2$, the Ising model, one computes   
$c(M_{\b})=(\tanh\b)^2$ which yields the correct reconstruction threshold. 

Theorem \ref{Theorem} is a generalization of the main result given in our paper 
\cite{FoKu08} for the specific case of 
the Potts model. 
That paper also contains comparisons of numerical values to the (presumed) 
exact transition values. Our discussion of the cases of $q=3,4,5$   
shows closeness up to a few percent, and for $q=3,4$ and small $d$  
these are the best rigorous bounds as of today.  
To see this connection between the present paper and   \cite{FoKu08} we rewrite 
$c(M_{\b})=\frac{e^{2\b}-1}{e^{2\b}+ q -1}\bar c(\b,q)$ and note that 
the main Theorem of  \cite{FoKu08} was formulated in terms of the quantity 
\begin{equation}\label{cbar}
\begin{split}
\bar c(\b,q)=
\sup_{p\in P }
\frac{\sum_{i=1}^q ( q p(i) -1)\log (1+(e^{2 \b}-1) p(i))}{\sum_{i=1}^q ( q p(i) -1)\log q p(i)}\cr
\end{split}.
\end{equation}

{\bf Numerical Example 2 (Non-symmetric Ising model.)}  
Consider the following transition matrix
\begin{equation}\label{1a}
  M=\left( \begin{matrix} 
      1-\d_1 & \d_1 \\
       1-\d_2& \d_2 \\
   \end{matrix}
  \right)\:\:\:\ \mbox{with}\:\:\:\d_1,\d_2\in{[0,1]}.\end{equation}
  The chain is not symmetric when $1-\d_1\neq\d_2$. 
  Let us focus on regular trees. Mossel and Peres in \cite{MoPe03} prove that, 
on a regular tree of degree $d$ the reconstruction problem defined by the matrix \eqref{1a} is unsolvable when 
  \begin{equation}\label{2a}
  d\frac{(\d_2-\d_1)^2}{\min\{\d_1+\d_2, 2-\d_1-\d_2\}}\leq 1,
  \end{equation}
{ while Martin in \cite{Martin03} gives the following condition for non-reconstructibility 
\begin{equation}\label{11a}
d\left(\sqrt{(1-\d_1)\d_2}-\sqrt{(1-\d_2)\d_1}\right)^2\leq 1.
\end{equation}}  

By the Kesten-Stigum bound it is known that there is reconstruction when $d(\d_2-\d_1)^2>1$.
When $\d_1+\d_2=1$, the matrix $M$ is symmetric and the Kesten--Stigum bound is sharp. Recently, Borgs, Chayes, Mossel and Roch 
in \cite{Bo06} have shown with an elegant proof that the Kesten--Stigum threshold is tight for roughly
 symmetric binary channels; i.e. when $|1-(\d_1+\d_2)|<\d$, for some $\d$ small.  Even if the threshold we give is very close to 
 Kesten--Stigum bound when the chain has a small asymmetry, by now, we are not able to recover this sharp estimate with our method. 
 For large asymmetry the Kesten--Stigum bound has been proved to not hold: Mossel proves as Theorem 1 
in \cite{Mo01} that, for any $\l>\frac{1}{d}$ there exists a $\d(\l)$ such that 
there is reconstruction for $\d_1,\d_2=\l+\d_1$ when $\d_1<\d$. 
 { On a Cayley tree with coordination number $d$, non-reconstruction for the Markov chain \eqref{1a} with $\d_2=0$ (or $1-\d_1=0$)
 is equivalent to the extremality of the Gibbs measure
for the hard-core model with activity $\frac{\d_1}{1-\d_1}\left(\frac{1}{1-\d_1}\right)^d$. Restricted to this specific case,
Martin proves a better condition than the one obtained taking $\d_2=0$ both in \eqref{11a} and in \eqref{C1}. }\\

Our entropy method  { provides a better bound than \eqref{11a}} and  considerably improves \eqref{2a} for the values of $\d_1$ and $\d_2$ giving a strongly asymmetric chain.\\
 
 A computation gives
 \begin{equation}
 c(M)=\sup_p\frac{(\d_2-\d_1)\log\left(\frac{(1-\d_2)+p(\d_2-\d_1)}{\d_2-p(\d_2-\d_1)}\frac{\d_1}{1-\d_2}\right)}{\log\left(\frac{p}{1-p}\frac{\d_1}{1-\d_2}\right)}.
 \end{equation}
 It is quite simple to compute numerically the constant $c(M)$; the numerical outputs and the comparisons with \eqref{2a}, {\red \eqref{11a}} and the Kesten-Stigum bound
 are in table \ref{tab:boundas}. 
 For the particular  
 pairs of values of $(\d_1,\d_2)$ we checked,  
 the Kesten-Stigum upper bound on the non-reconstruction thresholds for asymmetric chains are quite close to our lower bounds. 
 \begin{table}[h]
\centering
   \begin{tabular}{|c|c|c|c|c|} 
   \hline 
   $\d_1=0.3$&KS&FK&M&MP \\
   &\scriptsize{Kesten-Stigum}&\scriptsize{Formentin-K\"ulske}&\scriptsize{{ Martin}}&\scriptsize{Mossel-Peres}\\
   \hline
   $\d_2=0.1$&0.04&0.0579&{0.065}& 0.1\\
   \hline
   $\d_2=0.2$&0.01&0.0125&{ 0.0134}&0.02\\
    \hline 
   $\d_2=0.4$&0.01&0.0107&{0.0110}&0.0143\\
   \hline
   $\d_2=0.5$&0.04&0.0413&{ 0.0417}&0.05\\
   \hline
   $\d_2=0.6$&0.09&0.0907&{ 0.0910}&0.1\\
   \hline
   $\d_2=0.7$&0.16&0.16&0.16&0.16\\
   \hline
   $\d_2=0.8$&0.25&0.2525&{ 0.2534}&0.28\\
   \hline
   $\d_2=0.9$&0.36&0.3787&{ 0.3850}&0.45\\
   \hline
   \end{tabular}
   \bigskip\\
   \caption{For $\d_1=0.3$, the Kesten-Stigum upper bound on the non-reconstruction thresholds for asymmetric chains are very close to ours. 
   }
   \label{tab:boundas}
   \end{table}
   
%

%




\bigskip


\section{Proof}\label{s3}

We denote by $T^N$ the tree rooted at $0$ of depth $N$. 
The notation $T^N_v$ indicates the sub-tree of $T^N$ rooted at $v$ obtained from ``looking to the outside" on the tree $T^N$.  
We denote by $\mathbb{P}^{N}_v$ the measure on 
$T^N_v$ with free boundary conditions, or, equivalently the Markov chain obtained 
from broadcasting on the subtree with the root $v$ with the same transition kernel, starting 
in $\a$. 
We denote by $\mathbb{P}^{N,\xi}_v$ the correponding measure on 
$T^N_v$ with boundary condition on 
$\partial T^N_v$ given by $\xi=(\xi_i)_{i\in \partial T^N_v}$. Obviously it is obtained by conditioning 
the free boundary condition measure $\mathbb{P}^{N,\xi}_v$ to take the value $\x$ on the boundary. 

We write 
\begin{equation}\label{2.2}
\begin{split}
\pi^{N}_v=\pi^{N,\x}_v=\left(\mathbb{P}^{N,\xi}(\eta(v)=s)\right)_{s=1,\dots,q}\:.
\end{split}
\end{equation}

To control a recursion for these quantities along the tree we 
find it useful to make explicit the following notion. 

\begin{defn} We call a real-valued function $\LL$ on $P$ a linear stochastic  
Lyapunov function with center $p^*$ if there is a constant $c$ such that 
\begin{itemize}
\item $\LL(p)\geq 0\,\, \forall p \in P$ with equality if and only if $p=p^*$; 
\item $\E \LL(\pi^N_v) \leq c \sum_{w: v\mapsto w} \E \LL(\pi^{N}_w)$.
\end{itemize}
\end{defn}

\begin{prop}\label{proplja} Consider a tree-indexed Markov chain $\P$, with transition kernel 
$M(i,j)$ and invariant measure $\a(i)$. 

Then the function 
\begin{equation}\label{defL}
L(p)=S(p|\a)+S(\a|p)=
{ \sum_{i=1}^q \left(p(i)-\a(i)\right)\log\frac{p(i)}{\a(i)}}
\end{equation} 
is a linear stochastic  Lyapunov function with center $\a$ 
w.r.t. the measure $\P$ for the constant \eqref{C1}.
\end{prop}

Proposition \ref{proplja} immediately follows from the following 
invariance property of the recursion which is the main result of our paper.  

\begin{prop} \label{schoenerecursion} Main Recursion Formula for expected symmetrized entropy.
\begin{equation}
\begin{split}
\int \P(d\xi)L(\pi^{N,\xi}_v)= \sum_{w:v\rightarrow w}\int \P(d\xi) L(\pi^{N,\xi}_w M^\text{rev}).
\end{split}
\end{equation}
\end{prop} 
\bigskip
\bigskip
{\bf Warning: } Pointwise, that is for fixed boundary condition, things fail and one has  
\begin{equation}
\begin{split}
L(\pi^{N,\x}_v)\neq \sum_{w:v\rightarrow w}L(\pi^{N,\x}_w M^\text{rev})
\end{split}
\end{equation}
in general. In this sense the proposition should be seen as an invariance property 
which limits the possible behavior of the recursion. 

\bigskip

{\bf Proof of Proposition \ref{schoenerecursion}.}
We need the 
{\em measure on boundary configurations} at distance $N$ from the root on the tree 
emerging from $v$ which is obtained by 
conditioning the spin in the site $v$ to take the value to be $j$, namely \bigskip 
\begin{equation}
\begin{split}
Q^{N,j}_v(\xi):=\mathbb{P}^N_v(\s:\s_{|\partial T_v^N}=\xi|\:\s(v)=j) .
\end{split}
\end{equation}
Then the double expected value w.r.t. to the a priori measure $\a$ 
between boundary relative entropies can be written 
as an expected value w.r.t. $\P$ over boundary conditions w.r.t. 
to the open b.c. measure 
of the symmetrized entropy between the distributions at $v$ and 
$\a$ in the following form. \bigskip

\begin{lem} \label{lemma1} 
\begin{equation}
\begin{split}
&\int \P(d\xi)\underbrace{L(\pi^{N,\x}_v)}_{\text{symmetric entropy at } v}
=\int\a(d x_1)\int \a(d x_2) \underbrace{S(Q_v^{N,x_2}| Q_v^{N,x_1})}_{\text{boundary entropy}}.\cr 
\end{split}
\end{equation}
\end{lem}
\bigskip 

{\bf Proof of Lemma \ref{lemma1}: } 
In the first step we express the relative entropy as an expected value 
\begin{equation}
\begin{split}
S(Q_v^{N,x_2}| Q_v^{N,x_1})=\int \P(d\xi)\frac{d\pi^N_v}{d\a}(x_2)
\Bigl(\log\frac{d\pi^N_v}{d\a}(x_2)
-\log\frac{d\pi^N_v}{d\a}(x_1)\Bigr)\cr
\end{split}.
\end{equation}
Here we have used that, with obvious notations,   
\begin{equation}\label{tenn}
\begin{split}
\frac{dQ^{N,x_2}_v}{d\mathbb{P}_v^N}(\xi)=
\frac{\P_v(\s(v)=x_2,\x)}{\P_v(\s(v)=x_2)\P_v(\x)}
=\frac{d\pi^N_v}{d\a}(x_2). 
\end{split}
\end{equation}
Further we have used that
\begin{equation}\label{elf}
\begin{split}
&\log \frac{dQ^{N,x_2}_v}{dQ^{N,x_1}_v}
=\log\frac{d\pi^N_v}{d\a}(x_2)
-\log\frac{d\pi^N_v}{d\a}(x_1)\cr 
\end{split},
\end{equation}
for $x_1,x_2\in \{1,\dots,q\}$. This gives 
\begin{equation}
\begin{split}
&\int\a(d x_1)\int \a(d x_2) S(Q_v^{N,x_2}| Q_v^{N,x_1})\cr
&=\int \P(d\xi)\int\a(d x_2)\frac{d\pi^N_v}{d\a}(x_2)\log\frac{d\pi^N_v}{d\a}(x_2)\cr
&-\int \P(d\xi)\int\a(d x_1)\underbrace{\int\a(d x_2)\frac{d\pi^N_v}{d\a}(x_2)}_{1}\log\frac{d\pi^N_v}{d\a}(x_1)\cr
&=\int \P(d\xi)S(\pi^{N,\x}_v|\a)+\int \P(d\xi)S(\a | \pi^{N,\xi}_v)\cr
\end{split}
\end{equation}
and finishes the proof of Lemma \ref{lemma1}. 
$\Cox$

Let us continue with the proof of the Main Recursion Formula. 
We need two more ingredients formulated in the next two lemmas.  
The first gives the recursion of the probability vectors $\pi_{v}^N$ 
in terms of the values $\pi_{w}^N$ of their children $w$, which is valid 
for any fixed choice of the boundary condition $\x$. 


\begin{lem} \label{lemma2}  Deterministic recursion. 
\begin{equation}
\label{iteration}
\pi_{v}^N(j)=\frac{
\a(j)\prod_{w:v\rightarrow w}\sum_{i}\frac{M(j,i)}{\a(i)}\pi_{w}^N(i)}
{\sum_k\a(k)\prod_{w:v\rightarrow w}\sum_{i}\frac{M(k,i)}{\a(i)}\pi_{w}^N(i)},
\end{equation}
or, equivalently: 
for all pairs of values  $ j, k $ we have 
\begin{equation}
\label{iteration2}
\log\frac{d\pi^N_v}{d\a}(j)
-\log\frac{d\pi^N_v}{d\a}(k)
=\sum_{w:v\rightarrow w}
\log\frac{\sum_{i}\frac{M(j,i)}{\a(i)}\pi_{w}^N(i)
}
{\sum_{i}\frac{M(k,i)}{\a(i)}\pi_{w}^N(i)}\:.
\end{equation}
\end{lem}

The proof of this Lemma follows from an elementary computation 
with conditional probabilities and will be omitted here.  

We also need to take into account the {\em forward propagation} of the distribution 
of boundary conditions from the parents to the children, formulated in the next lemma. 

\begin{lem} \label{lemma3}  Propagation of the boundary measure.  
\begin{equation}\label{eq21}
\begin{split}
Q^{N,j}_v&=\prod_{w: v\rightarrow w} \sum_i M(j ,i ) Q^{N,i}_{w} \cr 
\end{split}.
\end{equation}
\end{lem} 

This statement follows from the definition of the model. 
Now we are ready to head for the Main Recursion Formula.

We use 
the second form of the statement of the deterministic recursion 
Lemma \ref{lemma2} equation \eqref{elf}
to write the boundary entropy in the form
\begin{equation}
\begin{split}
&S(Q_v^{N,j}| Q_v^{N,k})
=Q^{N,j}_v \sum_{w:v\rightarrow w}\log\frac{\sum_{i}\frac{M(j,i)}{\a(i)}\pi_{w}^N(i)
}
{\sum_{i}\frac{M(k,i)}{\a(i)}\pi_{w}^N(i)}\cr
\end{split}\:.
\end{equation}

Next, substituting the Propagation-of-the-boundary-measure-Lemma \ref{lemma3}  
and \eqref{tenn} we write 
\begin{equation}
\begin{split}\label{sia}
&S(Q_v^{N,j}| Q_v^{N,k})
=Q^{N,j}_v \sum_{w:v\rightarrow w}\log\frac{\sum_{i}\frac{M(j,i)}{\a(i)}\pi_{w}^N(i)
}
{\sum_{i}\frac{M(k,i)}{\a(i)}\pi_{w}^N(i)}\cr
&= \sum_{w:v\rightarrow w}\sum_l M(j ,l ) Q^{N,l}_{w} 
\log\frac{\sum_{i}\frac{M(j,i)}{\a(i)}\pi_{w}^N(i)
}
{\sum_{i}\frac{M(k,i)}{\a(i)}\pi_{w}^N(i)}
\cr
&= \sum_{w:v\rightarrow w} \int d\P(\xi)\sum_l  M(j ,l ) 
\frac{\pi^{N}_{w}(l)}{\a(l)}\log\frac{\sum_{i}\frac{M(j,i)}{\a(i)}\pi_{w}^N(i)
}
{\underbrace{\sum_{i}\frac{M(k,i)}{\a(i)}\pi_{w}^N(i)}_{\frac{\pi^N_w M^{\text{rev}}(k)}{\a(k)} }}
\cr
&= \sum_{w:v\rightarrow w} \int d\P(\xi)
\frac{\pi^N_w M^{\text{rev}}(j)}{\a(j)}
\log\frac{\frac{\pi^N_w M^{\text{rev}}(j)}{\a(j)}}
{\frac{\pi^N_w M^{\text{rev}}(k)}{\a(k)}},
\cr
\end{split}
\end{equation}
using in the last step the definition of the reversed Markov chain. 
Finally applying the sum $\sum_{j,k}\a(j)\a(k) \cdots $ to both 
sides of \eqref{sia} we get the Main Recursion Formula. To see this,  
note that the l.h.s. of \eqref{sia} together with this sum 
becomes the r.h.s. of the equation in Lemma \ref{lemma1}. 
For the r.h.s. of \eqref{sia} we note that 
\begin{equation}
\begin{split}\label{siam}
\sum_{j,k}\a(j)\a(k) 
\frac{\pi^N_w M^{\text{rev}}(j)}{\a(j)}
\log\frac{\frac{\pi^N_w M^{\text{rev}}(j)}{\a(j)}}
{\frac{\pi^N_w M^{\text{rev}}(k)}{\a(k)}}=L(\pi^N_w M^{\text{rev}})
\cr
\end{split}.
\end{equation}
This finishes the proof of the Main Recursion Formula Proposition \ref{schoenerecursion}. 
$\Cox$
\bigskip

Finally,  Theorem 
\ref{Theorem} follows from Proposition 
\ref{proplja} with the aid of the Wald equality with respect to the expectation 
over Galton-Watson trees since the contraction of the recursion 
and the Lyapunov function properties yield 
\begin{equation}\label{2.2666}
\begin{split}
\lim_{N\uparrow \infty}
\P\Bigl(\xi :  \Bigl |\pi^{N,\xi}(s)- \a(s) \Bigr | \geq \e\Bigr) \rightarrow 0 
\end{split},
\end{equation}
for all $s$, for all $\e>0$,  and this implies the extremality of the measure $\P$. 
This ends the proof of Theorem 
\ref{Theorem}.$\Cox$ 

\bigskip
\bigskip  

\textbf{Acknowledgements: }\\
The authors thank Aernout van Enter for interesting discussions and a critical reading of the manuscript.


\begin{thebibliography}{}
\bibitem{Io96a} D. Ioffe, {\em A note on the extremality of the disordered state for the Ising 
model on the Bethe lattice}, Lett. Math. Phys. \textbf{37}, 137-143 (1996).
\end{thebibliography}
\end{document}